\documentclass[a4paper,11pt]{article}
\usepackage{amsmath}
\usepackage{amssymb}
\usepackage{amsthm}
\usepackage{graphicx}
\usepackage{amscd}
\usepackage{epic, eepic}
\usepackage{url}
\usepackage{color}
\usepackage[utf8]{inputenc} 
\usepackage{comment}
\usepackage{enumerate}   
\usepackage{graphicx}
\usepackage{epstopdf}
\usepackage{enumitem}
\usepackage{tikz}
\usepackage[top=24mm, bottom=25mm, left=25mm, right=25mm]{geometry}
\usepackage[bookmarks=false,hidelinks]{hyperref}
\usepackage{mathrsfs}
\usetikzlibrary{arrows}
\definecolor{wrwrwr}{rgb}{0.3803921568627451,0.3803921568627451,0.3803921568627451}
\newcommand{\floor}[1]{\left\lfloor{#1}\right\rfloor}

\newtheorem{theorem}{Theorem}

\newtheorem{corollary}{Corollary}

\newtheorem{remark}{Remark}
\newtheorem{definition}{Definition}

\def\F{\mathcal F}
\newcommand{\ex}{{\rm  ex}}

\usepackage{authblk}

\title{Extremal problems for a matching and any other graph}
\author{Xiutao Zhu}
\author[]{Yaojun Chen}
\date{}

\affil[]{Department of Mathematics, Nanjing University, Nanjing 210093, P.R. CHINA.}

\begin{document}

\maketitle
\begin{abstract}
For a family of graphs $\F$, a graph is called $\F$-free if it does not contain any member of $\F$ as a subgraph. The generalized Tur\'an number $\ex(n,K_r,\F)$ is the maximum number of $K_r$ in an $n$-vertex $\F$-free graph and $\ex(n,K_2,\F)=\ex(n,\F)$, i.e., the classical Tur\'an number. Let $M_{s+1}$ be a matching on $s+1$ edges and $F$ be any graph. In this paper, we determine $\ex(n,K_r, \{M_{s+1},F\})$ apart from a constant additive term and also give a condition when the error constant term can be determined. In particular, we give the exact value of $\ex(n,\{M_{s+1},F\})$ for $F$ being any non-bipartite graph or some bipartite graphs.
Furthermore, we determine $\ex(n,K_r,\{M_{s+1},F\})$ when $F$ is color critical with $\chi(F)\ge \max\{r+1,4\}$.
 These extend the results in \cite{2023Alon,2023Gerbner,2023Ma}.

\end{abstract}

\section{Introduction}

 In this paper, let $K_r, K_{s,t}$ and $S_k$ denote the complete graph on $r$ vertices, complete bipartite graph with two parts of size $s$ and $t$, a star on $k$ vertices, respectively. Let $M_s$ denote a matching on $s$ edges. We use $|G|$ to denote the number of vertices of $G$. For a family of graphs $\F$, a graph is called $\F$-free if it does not contain any member of $\F$ as a subgraph.
 Let $G(n,s,H)$ denote the graph obtained from the complete bipartite graph $K_{s,n-s}$ by embedding a maximum $H$-free graph into the part of size $s$. For a subset $U\in V(G)$, we use $G[U]$ and $G-U$ to denote the subgraph induced by $U$ and $V(G)-U$, respectively.

 
  For an integer $r$ and a family $\F$, the generalized Tur\'an number $\ex(n,K_r,\F)$ is the maximum number of copies of $K_r$ in an $n$-vertex $\F$-free graph. Note that $\ex(n,K_2,\F)=\ex(n,\F)$, i.e., the classical Tur\'an number.  
The generalized Tur\'an number was firstly proposed by Alon and Shikhelman \cite{ALON2016146} in 2016. It has received a lot of attention in the past few years. Many classical results on Tur\'an problem have been extended to generalized Tur\'an number and some other interesting problem are studied too, see \cite{star,Beka,Turangood,EVEncycle,Andrzej,Hatami,Luo,MR4014346,Wang,linearforest,Dirac}.

In this paper, we mainly focus on the Tur\'an problem concerning matching. The first result in this issue dates back to Erd\H{o}s and Gallai \cite{erdHos1959maximal}, they proved 
$$\ex(n,M_{s+1})=\max\{e(G(n,s,K_{s+1})),e(K_{2s+1})\},$$ and determined the extremal graphs. This result was extended to generalized Tur\'an number $\ex(n,K_r,M_{s+1})$ by Wang \cite{Wang}. Beyond that, Chv\'atal and  Hanson \cite{V1976Degrees}, and independently  by Balachandran and Khare \cite{2009Graphs}  using different method,  determined the value of $\ex(n,\{M_{s+1}, S_{k+1}\})$ (The case for $s=k$ was proved early by Abbott,  Hanson and  Sauer \cite{1972Intersection}).
Recently, Alon and Frankl \cite{2023Alon} suggested to study $\ex(n,\{M_{s+1},F\})$ for any $F$. If there is an edge $e$ in $F$ such that $\chi(F-e)<\chi(F)$, then we call $F$ a color critical graph. They obtained the following results.
\begin{theorem}(Alon and Frankl \cite{2023Alon})\label{alon}
\begin{enumerate}
    \item For all $n\ge 2s+1$, $\ex(n,\{M_{s+1},K_{k+1}\})=\max\{e(K_{2s+1}),e(G(n,s,K_k))\}$.
    \item Let $F$ be a color critical graph with $\chi(F)=k+1\ge 3$.  When $s$ is large and $n\gg s$, 
    $$\ex(n,\{M_{s+1},F\})=e(G(n,s,K_k)).$$
\end{enumerate}    
\end{theorem}

Follow these results, Gerbner \cite{2023Gerbner} constructed some possible lower bounds of $\ex(n,\{M_{s+1},F\})$ and determined $\ex(n,\{M_{s+1},F\})$ apart from a constant additive term for some special bipartite graph $F$.  
\begin{theorem}(Gerbner\cite{2023Gerbner})\label{Gerbner}
Let $F$ be a bipartite graph and $p$ be the smallest size of a color class in any proper 2-coloring of $F$ with $p\le s$. Then $$\ex(n,\{M_{s+1},F\})=(p-1)n+O(1).$$ 
\end{theorem}

It appears likely that for $r \ge 3$, the function $\ex(n,K_r,\{M_{s+1},F\})$ behaves very differently from the classical Tur\'an number $\ex(n,\{M_{s+1},F\})$. The result of (1) in Theorem \ref{alon} is also extended to the generalized Tur\'an number $\ex(n,K_r,\{M_{s+1},K_{k+1}\})$ by Ma and Hou \cite{2023Ma}.  Let $\mathcal{N}_r(G)$ denote the number of copies of $K_r$ in $G$.
\begin{theorem}(Ma and Hou \cite{2023Ma})\label{Hou-ma}
For  $n\ge 2s+1$ and $k\ge r\ge 2$,
\[\ex(n,K_r, \{M_{s+1},K_{k+1}\})=\max\{\mathcal{N}_r(K_{2s+1}),\mathcal{N}_r(G(n,s,K_k))\}.\]
\end{theorem}
Furthermore, they also provided some possible lower bounds for $\ex(n,K_r,\{M_{s+1},F\})$ and  asked the exact value of $\ex(n,K_r,\{M_{s+1},F\})$ when $\chi(F)\ge 3$.

In this paper, we consider the generalized Tur\'an number about the matching and another graph. Before showing our results, we need some definitions.
A covering $S$ of $F$ is a subset of $V(F)$ such that $F-S$ is an empty graph, i.e., there is no edge in $F-S$. Let $F$ be a graph and $p$ be an integer, we define a family of subgraphs $\F[p]$ as follow,
\begin{definition}\label{def1}
 If $F$ has no covering of size at most $p$, then $\F[p]=\{K_{p+1}\}$. Otherwise $\F[p]=\{F[S]: S~\text{is a covering of F with}~|S|\le p \}.$   
\end{definition}
In addition to this, we call the covering $S$ an independent covering if $S$ is an independent set in $F$. We also need the definition about the size of the minimum independent covering.
\begin{definition}
If $F$ is bipartite, then $p(F)=\min\{|S|:S~\text{is an independent covering of F }\}.$ If $\chi(F)\ge 3$, then $p(F)=\infty$. 
\end{definition}

Note that, if $F$ is a bipartite graph, then $p(F)$ is exactly the smallest size of a color class in any proper 2-coloring, as we mentioned in Theorem \ref{Gerbner}. We determine $\ex(n,K_r,\{M_{s+1},F\})$ apart from a constant additive term and some exact values of $\ex(n,K_r,\{M_{s+1},F\})$ for special $F$.
\begin{theorem} \label{Thm3}
Let $F$ be a graph and $M_{s+1}$ be a matching.  Let $p<\min\{s+1,p(F)\}$ and $\ex(p,K_{r-1},\F[p])$ attains the maximum at $p=t$. Then,
\[\ex(n,K_r,\{M_{s+1},F\})=\ex(t,K_{r-1},\F[t])n+O(1).\]
Moreover, if $p(F)\ge s+1$ and $t=s$, then 
\[\ex(n,K_r,\{M_{s+1},F\})=\ex(s,K_{r-1},\F[s])(n-s)+\ex(s,K_r,\F[s]),\]
and $G(n,s,\F[s])$ is the unique extremal graph.
\end{theorem}

\begin{remark}
 The function $\ex(p,K_{r-1},\F[p])$ is not necessarily increasing on $p$. A simple example is $F=C_5$ and $r=3$.  Since  $C_5$ has no covering of size $2$ but has a covering of size $3$, then $\F[2]=\{K_3\}$ and $K_2\cup K_1\in \F[p]$ for $p\ge 3$.  We get $\ex(2,K_2,\F[2])=1$ but $\ex(p,K_2,\F[p])=0$ for $p\ge 3$. Moreover, one can check that for all other odd cycle $C_k$,  $\ex(p,K_{r-1},\F[p])$ is not  increasing either.    
\end{remark} 

\section{Some applications of Theorem \ref{Thm3}}
\begin{corollary}
Suppose $p(F)\ge s+1$ and  $n$ is large enough, 
$$\ex(n,\{M_{s+1},F\})=s(n-s)+\ex(s,\F[s]).$$
Moreover, $G(n,s,\F[s])$ is the unique extremal graph.
\end{corollary}
\noindent\emph{\textbf{Proof.}} When we consider the classical Tur\'an number, then $r=2$ in Theorem \ref{Thm3}. Note that $\ex(p,K_1,\F[p])=p$ as long as the independent set $I_p$ is not in $\F[p]$. However since $p(F)\ge s+1$, that is to say $F$ has no independent covering of size less than $s+1$, $I_p\notin \F[p]$ when $p<\min\{s+1,p(F)\}$. Thus  $\ex(p,K_1,\F[p])$ attains the maximum at $p=s$.

Therefore, by Theorem \ref{Thm3}, $\ex(n,\{M_{s+1},F\})=s(n-s)+\ex(s,\F[s])$ and $G(n,s,\F[s])$ is the unique extremal graph. This corollary extends Theorem \ref{alon} and determined the exact value for all non-bipartite graphs and the bipartite graphs with $p(F)\ge s+1$. $\hfill\blacksquare$

\begin{corollary}
Suppose $p(F)\le s$, then
$$\ex(n,\{M_{s+1},F\})=(p(F)-1)n+O(1).$$
\end{corollary}
\noindent\emph{\textbf{Proof.}} Analogously, $\ex(p,K_1,\F[p])$ attains the maximum at $p=p(F)$, then by Theorem \ref{Thm3}, we are done. $\hfill\blacksquare$

\vskip 3mm
Using Theorem \ref{Thm3}, we also determine the generalized Tur\'an number $\ex(n,K_r,\{M_{s+1},F\})$ when $F$ is color critical.
\begin{theorem}
Let $F$ be a color critical graph with $\chi(F)=k+1\ge \max\{r+1,4\}$. When $s\ge c(F,r)$ and $n\gg s$,
$$\ex(n,K_r,\{M_{s+1},F\})=\ex(s,K_{r-1},K_k)(n-s)+\ex(s,K_r,K_k)$$
and $G(n,s,K_k)$ is the unique extremal graph.
\end{theorem}
\noindent\emph{\textbf{Proof.}} Since $\chi(F)\ge 4$, we have $p(F)=\infty$. To use Theorem \ref{Thm3}, we need to study the property of $\ex(p,K_{r-1},\F[p])$ for $p<\min \{s+1,p(F)\}=s+1$.

Since $\chi(F)=k+1$, all graphs in $\F[p]$ have chromatic number at least $k$. If not, the chromatic number of $F$ would not exceed $k$ by the definition of $\F[p]$. Hence we have 
$$\ex(p,K_{r-1},\F[p])\ge \mathcal{N}_{r-1}(T_{k-1}(p)),$$
here $T_{k-1}(p)$ denotes the balanced complete ($k-1$)-partite graph on $p$ vertices(called Tur\'an graph).

On the other hand, since $F$ is color critical, we can find a $k+1$-coloring with the color class $V_1,V_2,\dots,V_{k+1}$ such that there is only one edge between $V_1,V_2$. Then if we delete the color class $V_{k+1}$, the resulting graph $F^-=F-V_{k+1}$ is still color critical and $\chi(F^-)=k\ge 3$. By the following theorem,
\begin{theorem}(Ma and Qiu \cite{MR4014346})
 Let $H$ be a color critical graph with $\chi(H)=k>m\ge 2$. Then when $n\ge c_0(H,m)$,  $\ex(n,K_m,H)=\mathcal{N}_m(T_{k-1}(n))$.  
\end{theorem}
  If we take $H=F^-$ and $m=r-1$ in the above theorem, then we know when $p\ge c_0(F,r)$, $\ex(p,K_{r-1},F^-)=\mathcal{N}_{r-1}(T_{k-1}(p))$. Note that $F^-\in \F[p]$, then 
\[
\mathcal{N}_{r-1}(T_{k-1}(p)) \le \ex(p,K_{r-1},\F[p])\le \ex(p,K_{r-1},F^-)=\mathcal{N}_{r-1}(T_{k-1}(p)). 
\]
That is to say, $\ex(p,K_{r-1},\F[p])= \mathcal{N}_{r-1}(T_{k-1}(p))$  is an increasing function when $s\ge p\ge c_0(F,r)$. For $p\le c_0(F,r)$, $\ex(p,K_{r-1},\F[p])$ does not exceed a  large constant $C$. Thus we can let $s$ be a large constant depending on $(F,r)$ so that $\ex(s,K_{r-1},\F[s])$ attains the maximum. 

Then by Theorem \ref{Thm3}, we know $G(n,s,\F[s])$ is the unique extremal graph. This extends the Theorems \ref{alon} and \ref{Hou-ma}. $\hfill\blacksquare$

\section{Proof of Theorem \ref{Thm3}}
In this section, we prove Theorem \ref{Thm3}. Let  $p<\min\{s+1,p(F)\}$ and $\ex(p,K_{r-1},\F[p])$ attains the maximum at $p=t$. Then $G(n,t,\F[t])$ is $\{M_{s+1},F\}$-free by Definition \ref{def1} and 
\[\begin{split}
 \ex(n,K_r,\{M_{s+1},F\})&\ge \ex(t,K_{r-1},\F[t])(n-t)+\ex(t,K_r,\F[t])\\
&= \ex(t,K_{r-1},\F[t])n+O(1).   
\end{split}\]
So the lower bound is done.

Next we prove the upper bound. Let $G$ be the extremal graph of $\ex(n,K_r,\{M_{s+1},F\})$. We need the following well-known theorem to discuss the structure of $G$.

\begin{theorem}(Tutte-Berge \cite{Berge})\label{Tutte}
The  graph $G$ is $M_{s+1}$-free if and only if there is a subset $B\subseteq V(G)$ such that for all components $G_1,\dots,G_m$ of $G-B$, they satisfy
\begin{align}\label{eq1}
 |B|+\sum_{i=1}^m \floor{\frac{|G_i|}{2}}\le s.  
\end{align}

\end{theorem}

Since $G$ is $M_{s+1}$-free, there is a set $B$ satisfying the inequality (\ref{eq1}) in the above theorem. Let $G_1,\dots,G_m$ be all components of $G-B$. Note that $s$ is a fixed constant, then most of these components are isolated vertices. Without loss of generality, we may assume  the components $G_1,\dots,G_\ell$ are not isolated vertices. 

Let $N_j$ denote the number of copies of $K_r$ which have $j$ vertices in $V(G)-B$ and $r-j$ vertices in $B$. Obviously, $N_0\le \binom{s}{r}$. For other $2\le j\le r$, by inequality (\ref{eq1}), we have $|B|+\sum_{i=1}^\ell|G_i|\le 3s$ and hence
\[\begin{split}
N_j\le \sum_{i=1}^\ell \mathcal{N}_j(G_i)\binom{s}{r-j} <\binom{3s}{r}=O(1), 
\end{split}\]
the second inequality holds since we can view $B\cup G_1\cup \dots \cup G_{\ell}$ as a big clique. This implies
\begin{align}
\mathcal{N}_r(G)=\sum_{j=0}^rN_j=N_1+O(1).    
\end{align}

Now we mainly deal with the term $N_1$. We divide $V(G)\setminus B$ into many subsets by the following way: let $U$ be a subset of $B$,  \[A_U=\{v\in V(G)\setminus B:N(v)\cap B=U\}.\]  
 Let $R=\{U: |A_U|\ge |F|\}$ and $Q=\{U: |A_U|<|F|\}$. Note that $|R|+|Q|= 2^{|B|}$ and we have

 \begin{align}
  N_1=\sum_{U\in R}\mathcal{N}_{r-1}(G[U])|A_U|+\sum_{U\in Q}\mathcal{N}_{r-1}(G[U])|A_U|.   
 \end{align}

For the set $U$ in $Q$, we have $\mathcal{N}_{r-1}(G[U])|A_U|<\binom{|B|}{r-1}|F|$ and hence
\begin{align}
\sum_{U\in Q}\mathcal{N}_{r-1}(G[U])|A_U|<2^{|B|}|F|\binom{|B|}{r-1}\le 2^s|F|\binom{s}{r-1}=O(1).
\end{align}
For the set $U$ in $R$, since $|A_U|\ge |F|$ and $G[U,A_U]$ is a complete bipartite graph, we can deduce that $G[U]$ is $\F[|U|]$-free by Definition \ref{def1} and $|U|<\min\{s+1,p(F)\}$. Hence $$\mathcal{N}_{r-1}(G[U])|A_U|\le \ex(|U|,K_{r-1},\F[|U|])|A_U|.$$ On the other hand, as we assumed, $t$ is the integer less than $\min\{s+1,p(F)\}$ such that $\ex(p, K_{r-1},\F[p])$ attains the maximum, then
\begin{align}
\begin{split}
 \sum_{U\in R}\mathcal{N}_{r-1}(G[U])|A_U|\le & \sum_{U\in R} \ex(|U|,K_{r-1},\F[|U|])|A_U|\\
 \le &\ex(t,K_{r-1},\F[t])\sum_{U\in R}|A_U|\\
 \le &\ex(t,K_{r-1},\F[t])(n-|B|).
\end{split}
\end{align}
Now combine the inequality (2)-(5), we know $\mathcal{N}_r(G)=\ex(t,K_{r-1},\F[t])n+O(1)$. We complete the proof of the first part in Theorem \ref{Thm3}.

Next we prove the second part, at this time $p(F)\ge s+1$ and $\ex(p,K_{r-1},\F[p])$ attains the maximum at $p=s$. Furthermore, by the inequality (2)-(5) and the lower bound, we have
\[\ex(s,K_{r-1},\F[s])(n-s)\le \mathcal{N}_r(G)\le  \sum_{U\in R} \ex(|U|,K_{r-1},\F[|U|])|A_U|+O(1).\]
Therefore, when $n$ is large,  $B\in R$ with $|B|=s$ and the corresponding $A_B$ satisfying $A_B=n-O(1)$.
Otherwise the right hand of the above would not exceed $\ex(s,K_{r-1},\F[s])(n-s)$. Hence we know $G[B]$ is $\F[s]$-free since $A_B$ is large. This also implies all components of $G-B$ are isolated vertices by Theorem \ref{Tutte}. For other vertices in $V(G)-B\cup A_B$, we can add all missing edges between them and $B$, this would not create any copy of $F$ since $G[B]$ is $\F[s]$-free, but the number of $K_r$ increases. Thus, we know $G=G(n,s,\F(s))$ and we are done. 
$\hfill\blacksquare$

\section{Classical Tur\'an number for balanced forests}

By the corollaries in Section 2, the only unsolved case for classical Tur\'an number is that $F$ is a bipartite graph with $p(F)\le s$. In this section we  deal with this case. A tree $T[A,B]$ is balanced if $|A|=|B|$. A  forest is balanced if each of its component is a balanced tree. The Tur\'an number of balanced forest was studied firstly by  Bushaw and Kettle \cite{kettle} if this forest contains at least two components. Here we study the Tur\'an problem combining the matching and a balanced forest and give a unified proof no matter the forest contains how many components. 
\begin{theorem}\label{thm2}
Let $F$ be a balanced forest with $v(F)=2p\le 2s$. If Erd\H{o}s-S\'os conjecture holds for each component of $F$ and $n$ is large, then 
\[\ex(n,\{F, M_{s+1}\})=(p-1)(n-p+1)+\ex(p-1,\F[p-1]).
\]
If $F$ has at least two components, then $G(n,p-1,\F[p-1])$ is the unique extremal graph. If $F$ is a tree, then $G(n-t(2p-1),p-1,\F[p-1])\cup tK_{2p-1}$ is the extremal graph, where $t\le \frac{s-p+1}{(p-1)}$. 
\end{theorem}

\begin{remark}
    By a result of Bushaw and Kettle(see Lemma 3.4 and Lemma 3.5 in  \cite{kettle}), it is easy to prove that  $\ex(p-1,\F[p-1])=\binom{p-1}{2}$ if $F$ contains a perfect matching and $\ex(p-1,\F[p-1])=0$, otherwise.
\end{remark}

\noindent\emph{\textbf{Proof.}} Let $F=T_1\cup \cdots \cup T_k $ be a balanced forest. The graph $G(n,p-1,\F[p-1])$ has matching number $p-1$ and it is also $F$-free by the definition of $\F[p-1]$. Besides this, if $F$ is a tree, then  $G(n-t(2p-1),p-1,\F[p-1])\cup tK_{2p-1}$ is $F$-free and has matching number at most $(p-1)+t(p-1)\le s$.
So the lower bound is done.

Next we prove the upper bound. Let $M=\{v_1u_1,\dots,v_tu_t\}$ be a maximum matching in $G$, $t\le s$.   This implies $V(G-M)$ is an independent set. We divide $V(G-M)$ into two subsets  $W$ and $W'$ such that 
\[W'=\{v\in V(G-M):d(v)\ge p\}~\text{and}~W=\{v\in V(G-M):d(v)\le p-1\}.\]

First we claim that $|W'|\le p\binom{2t}{p}$. Indeed, since there are at most $\binom{2t}{p}$ $p$-sets in $V(M)$ and if there are $p\binom{2t}{p}$ vertices in $W'$, then by the pigeonhole principle, there are $p$ vertices in $W'$ such that they have $p$ common neighbors in $V(M)$. Then we find a large complete bipartite graph, and hence a copy of $F$, a contradiction.

Next we assert that there are at least $2s\binom{2t}{p}$ vertices of degree $p-1$ in $W$. If not, then 
\[\begin{split}
  e(G)\le & \binom{2t}{2}+2t|W'|+ 2s\binom{2t}{p}(p-1)+(p-2)\left(n-2t-|W'|-2s\binom{2t}{p}\right) \\
  <&(p-1)(n-p+1).
\end{split}\]
The last inequality holds when $n$ is large, a contradiction. This also implies we can find $2s$ vertices of degree $p-1$ in $W$ such that they have $p$ common neighbors in $M$. Without loss of generality, let these $2p$ vertices be $\{x_1,\dots,x_{2s}\}$ and $U=\{v_1,\dots,v_{p-1}\}$ be the set of the common neighbors of them. 

On the other hand, for all other vertices in $W$ and the vertices in $V(M)-U$ whose degree is at most $p-1$, we can change their neighborhoods to $U$. This operation does not decrease the number of edges  and the resulting  graph is still $\{F,M_{s+1}\}$-free. Since if there is a copy of $F$ or $M_{s+1}$ in the resulting graph,  then we can use the vertices in $\{x_1,\dots,x_{2s}\}$  to replace the vertices in the copy of $F$ or $M_{s+1}$ which are incident with some new edges. That is, we can find a copy of $F$ or $M_{s+1}$ in $G$, a contradiction.

After the operation, let us redefine the set $W$ and $W'$, where $W$ denotes the set of all vertices of degree $p-1$ and have the neighborhood $U$, $W'$ denote the other vertices of degree at least $p$.
By above, $W'$ consists of vertices in the original $W'$ and some vertices in $V(M)$ whose degree is at least $p$. Thus $|B|\le K$ for some constant $K$. Furthermore, $W$ is independent with $|W|\ge n-K-2t=n-O(1)$. Since $G[U,W]$ is a large  complete bipartite graph, we have $G[U]$ is $\F[p-1]$-free by the definition. 

We also claim that there is no edge between $U$ and $W'$. If not, suppose $v_iw$ is an edge between $U$ and $W'$. Recall that $F=F[A,B]$ is a bipartite graph and $A,B$ are two color classes. There is a vertex, saying $x$, whose all neighbors except one are leaves. Suppose this vertex $x$ is in $B$ and let $y$ be the neighbor of $x$ which is not a leaf. Now we can embed the vertex $x$ into $w$, embed the vertex $y$ into $v_i$, embed the other neighbors of $x$ into the neighbor of $w$, embed the other vertices of $A$ into $U$ and the other vertices of $B$ into $W$. This can be done since $|U|=p-1$ and $d(w)\ge p$. 
Finally, we find a copy of $F$, a contradiction.

Therefore, $W'$ induces some connected components of $G$.    Furthermore, $G[W']$ is $T_1$-free. Otherwise, a copy of $T_1$ in $W'$ together with a copy of $T_2\cup \cdots \cup T_k$  in $G[U,W]$ would construct a copy of $F$. Thus, we have
\[\begin{split}
 e(G)\le& \ex(p-1,\F[p-1])+(p-1)(n-p+1-|W'|) +\ex(|W'|,T_1)\\
 \le & \ex(p-1,\F[p-1])+(p-1)(n-p+1)-|W'|(p-1)+\frac{v(T_1)-2}{2}|W'|\\
 \le & \ex(p-1,\F[p-1])+(p-1)(n-p+1).
\end{split}\]
The last inequality holds under the assumption of Erd\H{o}s-S\'os conjecture. From the above inequalities, if $F$ is a real forest, then $(p-1)>(v(T_1)-2)/2$. So if  the equality holds, then $W'=\emptyset$ and $G=G(n,p-1,\F[p-1])$. If $F$ is a tree, then the equality holds if and only if $W'$ induces some disjoint cliques $K_{2p-1}$. But since $G$ is $M_{s+1}$-free, $W'$ induces  at most $\frac{s-p+1}{p-1}$ copies of $K_{2p-1}$. That is $G=G(n-(p-1)-t(2p-1),p-1,\F[p-1])\cup tK_{2p-1}$ with $t\le \frac{s-p+1}{p-1}$. The proof is completed.$\hfill\blacksquare$

\vskip 3mm

\section{Acknowledgements}
This Research was supported by NSFC under grant numbers 12161141003 and 11931006.

\bibliography{Referances.bib}

\end{document}